\documentclass[12pt,headrule]{article}

\usepackage{graphicx}

\def\re#1{\par\hangindent\parindent\indent\llap{#1\enspace}\ignorespaces}

\def\c{\centerline}\def\no{\noindent}

\begin{document}

\c{\large\bf A Multi-Space Model for}\vskip 3mm

\c{\large\bf Chinese Bids Evaluation with Analyzing}

\vskip 4mm \c{Linfan Mao}{\it

\vskip 3mm \c{\scriptsize Chinese Academy of Mathematics and
System Sciences, Beijing 100080, P.R.China}

\c{\scriptsize Guoxin Tendering Co.,LTD, Beijing 100044,
P.R.China}

\c{maolinfan@163.com}}

\vskip 5mm

\begin{minipage}{130mm}
\no{\bf Abstract.} {\small A tendering is a negotiating process
for a contract through by a tenderer issuing an invitation,
bidders submitting bidding documents and the tenderer accepting a
bidding by sending out a notification of award. As a useful way of
purchasing, there are many norms and rulers for it in the
purchasing guides of the World Bank, the Asian Development Bank,
$\cdots$, also in contract conditions of various consultant
associations. In China, there is a law and regulation system for
tendering and bidding. However, few works on the mathematical
model of a tendering and its evaluation can be found in
publication. The main purpose of this paper is to construct a
Smarandache multi-space model for a tendering, establish an
evaluation system for bidding based on those ideas in the
references $[7]$ and $[8]$ and analyze its solution by applying
the decision approach for multiple objectives and value
engineering. Open problems for pseudo-multi-spaces are also
presented in the final section.}

 \no{\bf Key Words:} \ {\small tendering, bidding, evaluation,
Smarandache multi-space, condition of successful bidding, decision
of multiple objectives, decision of simply objective,
pseudo-multiple evaluation, pseudo-multi-space. }

\no{\bf AMS(2000)}: {\small 90B50,90C35,90C90}
\end{minipage}

\vskip 8mm

\no{\bf \S $1.$ Introduction}\vskip 3mm

\no The tendering is an efficient way for purchasing in the market
economy. According to the {\it Contract Law of the People's
Republic of China} (Adopted at the second meeting of the Standing
Committee of the 9th National People's Congress on March 15,1999),
it is just a civil business through by a tenderer issuing a
tendering announcement or an invitation, bidders submitting
bidding documents compiled on the tendering document and the
tenderer accepting a bidding after evaluation by sending out a
notification of award. The process of this business forms a
negotiating process of a contract. In China, there is an interval
time for the acceptation of a bidding and becoming effective of
the contract, i.e., the bidding is accepted as the tenderer send
out the notification of award, but the contract become effective
only as the tenderer and the successful bidder both sign the
contract.

In the {\it Tendering and Bidding Law of the People's Republic of
China} (Adopted at the 11th meeting of the Standing Committee of
the 9th National People's Congress on August 30,1999), the
programming and liability or obligation of the tenderer, the
bidders, the bid evaluation committee and the government
administration are stipulated in detail step by step. According to
this law, the tenderer is on the side of raising and formulating
rulers for a tender project and the bidders are on the side of
response each ruler of the tender. Although the bid evaluation
committee is organized by the tenderer, its action is independent
on the tenderer. In tendering and bidding law and regulations of
China, it is said that any unit or person can not disturbs works
of the bid evaluation committee illegally. The action of them
should consistent with the tendering and bidding law of China and
they should place themselves under the supervision of the
government administration.

The role of each partner can be represented by a tetrahedron such
as those shown in Fig.$1$.

\includegraphics[bb=-5 5 400 170]{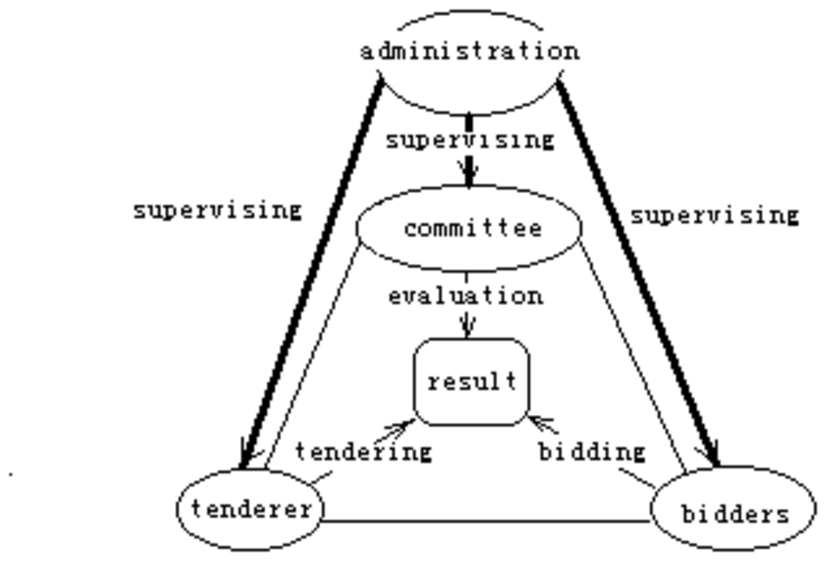}

\vskip 2mm

\c{\bf Fig.$1$}\vskip 2mm

The 41th item in the {\it Tendering and Bidding Law of the
People's Republic of China} provides conditions for a successful
bidder:

($1$) optimally responsive all of the comprehensive criterions in
the tendering document;

($2$) substantially responsive criterions in the tender document
with the lowest evaluated bidding price unless it is lower than
this bidder's cost.

The conditions ($1$) and ($2$) are often called the {\it
comprehensive evaluation method} and the {\it lowest evaluated
price method}. In the same time, these conditions also imply that
the tendering system in China is a multiple objective system, not
only evaluating in the price, but also in the equipments,
experiences, achievements, staff and the programme, etc.. However,
nearly all the encountered evaluation methods in China do not
apply the scientific decision of multiple objectives. In where,
the {\it comprehensive evaluation method} is simply replaced by
the {\it 100 marks} and the {\it lowest evaluated price method} by
the {\it lowest bidding price method}. Regardless of whether
different objectives being comparable, there also exist problems
for the ability of bidders and specialists in the bid evaluation
committee creating a false impression for the successful bidding
price or the successful bidder. The tendering and bidding is badly
in need of establishing a scientific evaluation system in
accordance with these laws and regulations in China. Based on the
reference $[7]$ for Smarandache multi-spaces and the mathematical
model for the tendering in $[8]$, the main purpose of this paper
is to establish a multi-space model for the tendering and a
scientific evaluation system for bids by applying the approach in
the multiple objectives and value engineering, which enables us to
find a scientific approach for tendering and its management in
practice. Some cases are also presented in this paper.

The terminology and notations are standard in this paper. For
terminology and notation not defined in this paper can be seen in
$[7]$ for multi-spaces, in $[1]-[3]$ and $[6]$ for programming,
decision and graphs and in $[8]$ for the tendering and bidding
laws and regulations in China.

\vskip 5mm

\no{\bf \S$2.$ A multi-space model for tendering }\vskip 3mm

\no Under an idea of anti-thought or paradox for mathematics :{\it
combining different fields into a unifying field}, Smarandache
introduced the conception of multi-spaces in 1969([9]-[12]),
including algebraic multi-spaces and multi-metric spaces. The
 contains the well-known Smarandache
geometries($[5]-[6]$), which can be used to {\it General
Relativity} and {\it Cosmological Physics}([7]). As an application
to {\it Social Sciences}, multi-spaces can be also used to
establish a mathematical model for tendering.

These algebraic multi-spaces are defined in the following
definition.

\vskip 4mm

\no{\bf Definition $2.1$} \ {\it An algebraic multi-space $\sum$
with multiple $m$ is a union of $m$ sets $A_1,A_2,\cdots,A_m$}

$$\sum = \bigcup\limits_{i=1}^mA_i,$$

\no{\it where $1\leq m< +\infty$ and there is an operation or
ruler $\circ_i$ on each set $A_i$ such that $(A_i£¬\circ_i)$ is an
algebraic system for any integer $i, 1\leq i\leq m$.}

\vskip 3mm

Notice that if $i\not=j, 1\leq i,j\leq m$, there must not be
$A_i\bigcap A_j=\emptyset$, which are just correspondent with the
characteristics of a tendering. Thereby, we can construct a
Smarandache multi-space model for a tendering as follows.

Assume there are $m$ evaluation items $A_1,A_2,\cdots, A_m$ for a
tendering $\widetilde{A}$ and there are $n_i$ evaluation indexes
$a_{i1},a_{i2},\cdots,a_{in_i}$ for each evaluation item $A_i,
1\leq i\leq m$. By applying mathematics, this tendering can be
represented by

$$\widetilde{A} = \bigcup\limits_{i=1}^mA_i,$$

\no where, for any integer $i, 1\leq i\leq m$,

$$(A_i,\circ_i)=\{a_{i1},a_{i2},\cdots,a_{in_i}|\circ_i\}$$

\no is an algebraic system. Notice that we do not define other
relations of the tendering $\widetilde{A}$ and evaluation indexes
$a_{ij}$ with $A_i,1\leq i\leq m$ unless $A_i\subseteq
\widetilde{A}$ and $a_{ij}\in A_i$ in this multi-space model.

Now assume there are $k,k\geq 3$ bidders $R_1,R_2,\cdots,R_k$ in
the tendering $\widetilde{A}$ and the bidding of bidder $R_j,1\leq
j\leq k$ is

\[
R_j(\widetilde{A})=R_j\left(
\begin{array}{l}
A_1 \\
A_2 \\
\cdots\\
A_m\\
\end{array}
\right) = \left(
\begin{array}{l}
R_j(A_1) \\
R_j(A_2) \\
\cdots\\
R_j(A_m)\\
\end{array}
\right)  \].\vskip 2mm

According to the successful bidding criterion in the {\it
Tendering and Bidding Law of the People's Republic of China} and
regulations, the bid evaluation committee needs to determine
indexes $i_1,i_2,\cdots,i_k$, where
$\{i_1,i_2,\cdots,i_k\}=\{1,2,\cdots,k\}$ such that there is an
ordered sequence

$$R_{i_1}(\widetilde{A})\succ R_{i_2}(\widetilde{A})\succ\cdots
\succ R_{i_k}(\widetilde{A})$$

\no for these bidding
$R_1(\widetilde{A}),R_2(\widetilde{A}),\cdots,R_k(\widetilde{A})$
of bidders $R_1,R_2,\cdots,R_k$. Here, these bidders
$R_{i_1},R_{i_2}$ and $R_{i_3}$ are {\it pre-successful bidders}
in succession determined by the bid evaluation committee in the
laws and regulations in China.

\vskip 4mm

\no{\bf Definition $2.2$} \ {\it An ordered sequence for elements
in the symmetry group $S_n$ on $\{1,2,\cdots,m\}$ is said an
alphabetical sequence if it is arranged by the following
criterions:

$(i)$ $(1,0\cdots,0)\succeq P$ for any permutation $P\in S_n$.

$(ii)$ if integers $s_1,s_2,\cdots,s_h\in\{1,2,\cdots,m\},1\leq h
< m$ and permutations $(s_1,s_2,$ $\cdots,s_h,t,\cdots),$
$(s_1,s_2,\cdots,s_h,l,\cdots)\in S_n$, then}

$$(s_1,s_2,\cdots,s_h,t,\cdots)\succ(s_1,s_2,\cdots,s_h,l,\cdots)$$

\no {\it if and only if $t< l$. Let $\{x_{\sigma_i}\}_1^n$ be a
sequence, where $\sigma_1\succ\sigma_2\succ\cdots\succ\sigma_n$
and $\sigma_i\in S_n$ for $1\leq i\leq n$, then the sequence
$\{x_{\sigma_i}\}_1^n$ is said an alphabetical sequence.

Now if $x_{\sigma}\succ x_{\tau}$, $x_{\sigma}$ is preferable than
$x_{\tau}$ in order. If $x_{\sigma}\succeq x_{\tau}$, then
$x_{\sigma}$ is preferable or equal with $x_{\tau}$ in order. If
$x_{\sigma}\succeq x_{\tau}$ and $x_{\tau}\succeq x_{\sigma}$,
then $x_{\sigma}$ is equal $x_{\tau}$ in order, denoted by
$x_{\sigma}\approx x_{\tau}$. }

\vskip 3mm

We get the following result for an evaluation of a tendering.

\vskip 4mm

\no{\bf Theorem $2.1$} \ {\it Let $O_1,O_2,O_3\cdots$ be ordered
sets. If $R_j(\widetilde{A})\in O_1\times O_2\times
O_3\times\cdots$ for any integer $j,1\leq j\leq k$, then there
exists an arrangement $i_1,i_2,\cdots,i_k$ for indexes
$1,2,\cdots,k$ such that}

$$R_{i_1}(\widetilde{A})\succeq R_{i_2}(\widetilde{A})
\succeq\cdots\succeq R_{i_k}(\widetilde{A}).$$

\vskip 3mm

{\it Proof} \ By the assumption, for any integer $j,1\leq j\leq
k,$

$$R_j({\widetilde{A}})\in O_1\times O_2\times O_3\times\cdots.$$

\no Whence, $R_j({\widetilde{A}})$ can be represented by

$$R_j({\widetilde{A}})=(x_{j1},x_{j2},x_{j3},\cdots),$$

\no where $x_{jt}\in O_t, t\geq 1$. Define a set

$$S_t=\{x_{jt}; 1\leq j\leq m\}.$$

\no Then the set $S_t\subseteq O_t$ is finite. Because the set
$O_t$ is an ordered set, so there exists an order for elements in
$S_t$. Not loss of generality, assume the order is

$$x_{1t}\succeq x_{2t}\succeq\cdots\succeq x_{mt},$$

\no for elements in $S_t$. Then we can apply the alphabetical
approach to $R_{i_1}(\widetilde{A}), R_{i_2}(\widetilde{A}),$
$\cdots, R_{i_k}(\widetilde{A})$ and get indexes
$i_1,i_2,\cdots,i_k$ such that

$$R_{i_1}(\widetilde{A})\succeq R_{i_2}(\widetilde{A})
\succeq\cdots\succeq R_{i_k}(\widetilde{A}). \ \ \natural$$

If we choose $O_i, i\geq 1$ to be an ordered function set in
Theorem $2.1$, particularly, let $O_1=\{f\}, f:A_i\rightarrow R,
1\leq i\leq m$ be a monotone function set and $O_t=\emptyset$ for
$t\geq 2$, then we get the next result.

\vskip 4mm

\no{\bf Theorem $2.2$} \ {\it Let $R_j:A_i\rightarrow R, 1\leq
i\leq m, 1\leq j\leq k$ be monotone functions. Then there exists
an arrangement $i_1,i_2,\cdots,i_k$ for indexes $1,2,\cdots,k$
such that}

$$R_{i_1}(\widetilde{A})\succeq R_{i_2}(\widetilde{A})
\succeq\cdots\succeq R_{i_k}(\widetilde{A}).$$

\vskip 3mm

We also get the following consequence for evaluation numbers by
Theorem $2.2$.

\vskip 4mm

\no{\bf Corollary $2.1$} \ {\it If $R_j(A_i)\in
[-\infty,+\infty]\times [-\infty,+\infty] \times
[-\infty,+\infty]\times\cdots$ for any integers $i,j,1\leq i\leq
m, 1\leq j\leq k$, then there exists an arrangement
$i_1,i_2,\cdots,i_k$ for indexes $1,2,\cdots,k$ such that}

$$R_{i_1}(\widetilde{A})\succeq R_{i_2}(\widetilde{A})
\succeq\cdots\succeq R_{i_k}(\widetilde{A}).$$

\vskip 3mm

Notice that in the above ordered sequence, if we arrange
$R_{i_s}\succ R_{i_l}$ or $R_{i_l}\succ R_{i_s}$ further in the
case of $R_{i_s}\approx R_{i_l}, s\not=l$, then we can get an
ordered sequence

$$R_{i_1}(\widetilde{A})\succ R_{i_2}(\widetilde{A})\succ\cdots
\succ R_{i_k}(\widetilde{A}),$$

\no and the pre-successful bidders accordance with the laws and
regulations in China.

\vskip 5mm

\no{\bf \S $3.$ A mathematical analog for bids evaluation}\vskip
3mm

\no For constructing an evaluation system of bids by the
multi-space of tendering, the following two problems should be
solved in the first.

\vskip 2mm

\no{\bf Problem $1$} \ {\it For any integers $i,j,1\leq i,j\leq
m$, how to determine $R_j(A_i)$ on account of the responsiveness
of a bidder $R_j$ on indexes $a_{i1},a_{i2},\cdots,a_{in_i}$?}

\no{\bf Problem $2$} \ {\it For any integer $j,1\leq j\leq m$, how
to determine $R_j(\widetilde{A})$ on account of the vector
$(R_j(A_1),R_j(A_2),\cdots,R_j(A_m))^t$?}

\vskip 2mm

Different approaches for solving Problems $1$ and $2$ enable us to
get different mathematical analogs for bids evaluation.

\vskip 3mm

\no{\bf $3.1.$ An approach of multiple objectives decision}

\vskip 3mm

\no This approach is originated at the assumption that
$R_j(A_1),R_j(A_2),\cdots,R_j(A_m), 1\leq j\leq m$ are independent
and can not compare under a unified value unit. The objectives of
tendering is multiple, not only in the price, but also in the
equipments, experiences, achievements, staff and the programme,
etc., which are also required by the 41th item in the {\it
Tendering and Bidding Law of the People's Republic of China}.

According to Theorems $2.1-2.2$ and their inference, we can
establish a programming for arranging the order of each evaluation
item $A_i,1\leq i\leq m$ and getting an ordered sequence of bids
$R_1({\widetilde{A}}),R_2(\widetilde{A}),$
$\cdots,R_k(\widetilde{A})$ of a tendering $\widetilde{A}=
\bigcup\limits_{i=1}^mA_i$, as follows:\vskip 2mm

{\bf STEP $1$}\ {\it determine the order of the evaluation items
$A_1,A_2,\cdots,A_m$. For example, for $m=5$, $A_1\succ A_2\approx
A_3\succ A_4\approx A_5$ is an order of the evaluation items
$A_1,A_2,A_3,A_4,A_5$.}

{\bf STEP $2$}\ {\it for two bids $R_{j_1}(A_i),R_{j_2}(A_i),
j_1\not=j_2, 1\leq i\leq m$, determine the condition for
$R_{j_1}(A_i)\approx A_{j_2}(A_2)$. For example, let $A_1$ be the
bidding price. Then $R_{j_1}(A_1)\approx R_{j_2}(A_1)$ providing
$|R_{j_1}(A)-R_{j_2}(A_1)|\leq 100$(10 thousand yuan).}

{\bf STEP $3$}\ {\it for any integer $i,1\leq i\leq m$, determine
the order of $R_1(A_i),R_2(A_i),$ $\cdots,R_k(A_i)$. For example,
arrange the order of bidding price from lower to higher and the
bidding programming dependent on the evaluation committee. }

{\bf STEP $4$}\ {\it alphabetically arrange
$R_1(\widetilde{A}),R_2(\widetilde{A}),\cdots,R_k(\widetilde{A})$,
which need an approach for arranging equal bids
$R_{j_1}(\widetilde{A})\approx R_{j_2}(\widetilde{A})$ in order.
For example, arrange them by the ruler of {\it lower price
preferable} and get an ordered sequence

$$R_{i_1}(\widetilde{A})\succ R_{i_2}(\widetilde{A})
\succ\cdots\succ R_{i_k}(\widetilde{A})$$

\no of these bids
$R_1(\widetilde{A}),R_2(\widetilde{A}),\cdots,R_k(\widetilde{A})$.}

Notice that we can also get an ordered sequence through by
defining the weight functions

$$\omega(\widetilde{A})=H(\omega(A_1),\omega(A_2),\cdots,\omega(A_m))$$

\no and

$$\omega(A_i)=F(\omega(a_{i1}),\omega(a_{i2}),\cdots,\omega(a_{in_i})).$$

\no For the weight function in detail, see the next section.

\vskip 4mm

\no{\bf Theorem $3.1$}\ {\it The ordered sequence of bids of a
tendering $\widetilde{A}$ can be gotten by the above programming.}

\vskip 3mm

{\it Proof} \ Assume there are $k$ bidders in this tendering. Then
we can alphabetically arrange these bids $R_{i_1}(\widetilde{A}),
R_{i_2}(\widetilde{A}),\cdots, R_{i_k}(\widetilde{A})$ and get

$$R_{i_1}(\widetilde{A})\succeq R_{i_2}(\widetilde{A})
\succeq\cdots\succeq R_{i_k}(\widetilde{A}).$$

Now applying the arranging approach in the case of
$R_{j_1}(\widetilde{A})\approx R_{j_2}(\widetilde{A})$, we finally
obtain an ordered sequence

$$R_{i_1}(\widetilde{A})\succ R_{i_2}(\widetilde{A})\succ
\cdots\succ R_{i_k}(\widetilde{A}). \ \ \natural$$

\vskip 3mm

\no{\bf Example $3.1$}\ {\it There are $3$ evaluation items in a
building construction tendering $\widetilde{A}$ with $A_1=$price,
$A_2$=programming and $A_3$=similar achievements in nearly $3$
years. The order of the evaluation items is $A_1\succ A_3\succ
A_2$ and $R_{j_1}(A_i)\approx R_{j_2}(A_i), 1\leq i\leq 3$
providing $|R_{j_1}(A_1)-R_{j_2}(A_1)|\leq 150$, $R_{j_1}(A_2)$
and $R_{j_2}(A_2)$ are in the same rank or the difference of
architectural area between $R_{j_1}(A_3)$ and $R_{j_2}(A_3)$ is
not more than $40000m^2$. For determining the order of bids for
each evaluation item, it applies the rulers that from the lower to
the higher for the price, from higher rank to a lower rank for the
programming by the bid evaluation committee and from great to
small amount for the similar achievements in nearly $3$ years and
arrange $R_{j_1}(\widetilde{A}), R_{j_2}(\widetilde{A})$, $1\leq
j_1,j_2\leq k=$bidders by the ruler of lower price first for two
equal bids in order $R_{j_1}(\widetilde{A})\approx
R_{j_2}(\widetilde{A})$.

There were $4$ bidders $R_1,R_2,R_3,R_4$ in this tendering. Their
bidding prices are in table $1$.

\vskip 3mm
\begin{center}
\begin{tabular}{|c|c|c|c|c|}\hline
\ bidder \ & $R_1$ &  $R_2$ & $R_3$ & $R_4$  \\ \hline
 \ $A_1$ \
 & $3526$ &  $3166$ & $3280$ & $3486$  \\ \hline
\end{tabular}
\end{center}
\vskip 2mm

\c{\rm table $1$}

Applying the arrangement ruler for $A_1$, the order for $R_2(A_1),
R_3(A_1), R_4(A_1),$ $R_1(A_1)$ is

$$R_2(A_1)\approx R_3(A_1)\succ R_4(A_1)\approx R_1(A_1).$$

The evaluation order for $A_2$ by the bid evaluation committee is
$R_3(A_2)\approx R_2(A_2)\succ R_1(A_2)\succ R_4(A_2)$. They also
found the bidding results for $A_3$ are in table $2$.

\vskip 3mm
\begin{center}
\begin{tabular}{|c|c|c|c|c|}\hline
\ bidder \ & $R_1$ &  $R_2$ & $R_3$ & $R_4$  \\ \hline
 \ $A_3(m^2)$ \
 & $250806$ &  $210208$ & $290108$ & $300105$  \\ \hline
\end{tabular}
\end{center}
\vskip 2mm \c{\rm table $2$}

Whence the order of $R_4(A_3), R_3(A_3), R_1(A_3), R_2(A_3)$ is

$$R_4(A_3)\approx R_3(A_3)\succ R_1(A_3)\approx R_2(A_3).$$

Therefore, the ordered sequence for these bids
$R_1(\widetilde{A}),R_2(\widetilde{A}),R_3(\widetilde{A})$ and
$R_4(\widetilde{A})$ is

$$R_3(\widetilde{A})\succ R_2(\widetilde{A})\succ
R_4(\widetilde{A})\succ R_1(\widetilde{A}).$$ }

Let the order of evaluation items be $A_1\succ A_2\succ\cdots\succ
A_m$. Then we can also get the ordered sequence of a tendering by
applying a graphic method. By the terminology in graph theory, to
arrange these bids of a tendering is equivalent to find a directed
path passing through all bidders $R_1,R_2,\cdots,R_k$ in a graph
$G[\widetilde{A}]$ defined in the next definition. Generally, the
graphic method is more convenience in the case of less bidders,
for instance $7$ bidders for a building construction tendering in
China.

\vskip 4mm

\no{\bf Definition $3.1$}\ {\it Let $R_1,R_2,\cdots,R_k$ be all
these $k$ bidders in a tendering
$\widetilde{A}=\bigcup\limits_{i=1}^mA_i$. Define a directed graph
$G[\widetilde{A}]=(V(G[\widetilde{A}]),E(G[\widetilde{A}]))$ as
follows.

$V(G[\widetilde{A}])=\{R_1,R_2,\cdots,R_k\}\times\{A_1,A_2,\cdots,A_m\},$

$E(G[\widetilde{A}])=E_1\bigcup E_2\bigcup E_3$.

\no Where $E_1$ consists of all these directed edges
$(R_{j_1}(A_i),R_{j_2}(A_i))$, $1\leq i\leq m,1\leq j_1,j_2\leq k$
and $R_{j_1}(A_i)\succ R_{j_2}(A_i)$ is an adjacent order. Notice
that if $R_s(A_i)\approx R_l(A_i)\succ R_j(A_i)$, then there are
$R_s(A_i)\succ R_j(A_i)$ and $ R_l(A_i)\succ R_j(A_i)$
simultaneously. $E_2$ consists of edges
$R_{j_1}(A_i)R_{j_2}(A_i),1\leq i\leq m,1\leq j_1,j_2\leq k$,
where $R_{j_1}(A_i)\approx R_{j_2}(A_i)$ and
$E_3=\{R_j(A_i)R_j(A_{i+1})|1\leq i\leq m-1, 1\leq j\leq k\}$.}

\vskip 3mm

For example, the graph $G[\widetilde{A}]$ for Example $3.1$ is
shown in Fig.$2$.

\vskip 3mm

\includegraphics[bb=50 5 400 230]{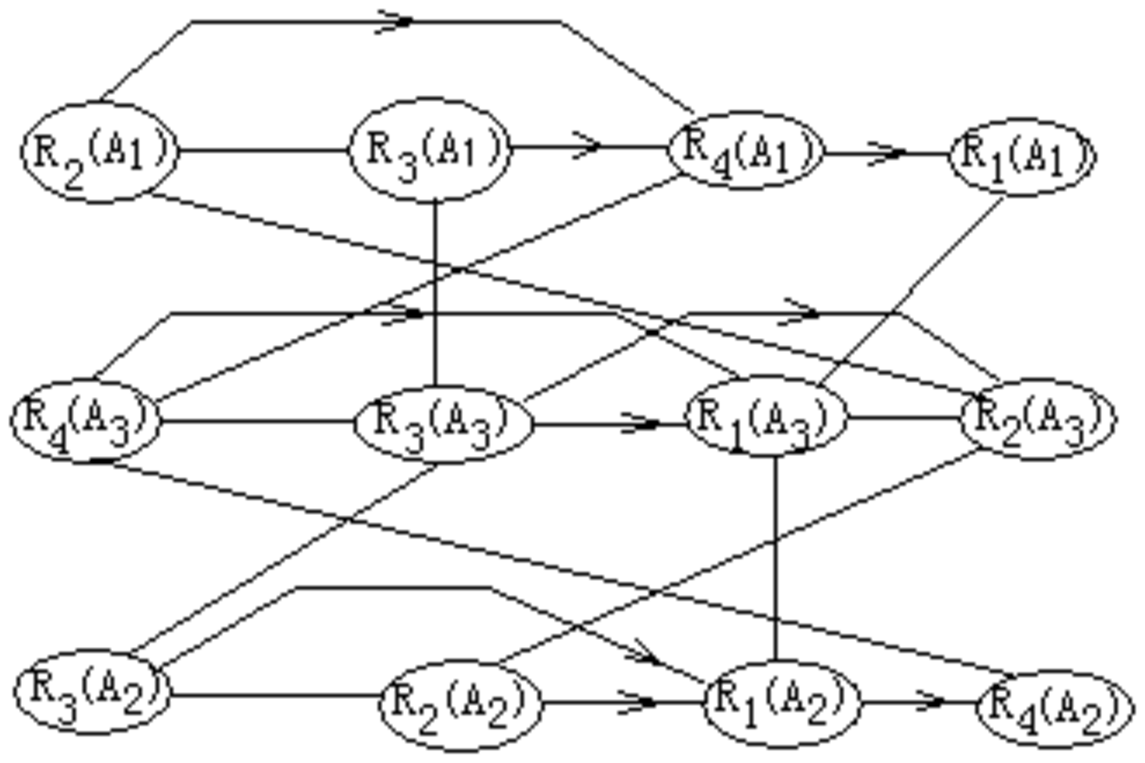}

\vskip 2mm

\c{Fig.$2$}

\vskip 2mm

Now we need to find a directed path passing through $R_1, R_2,
R_3,R_4$ with start vertex $R_2(A_1)$ or $R_3(A_1)$. By the ruler
in an alphabetical order, we should travel starting from the
vertex $R_3(A_1)$ passing through $A_2,A_3$ and then arriving at
$A_1$. Whence, we find a direct path correspondent with the
ordered sequence

$$R_3(\widetilde{A})\succ R_2(\widetilde{A})\succ
R_4(\widetilde{A})\succ R_1(\widetilde{A}).$$

\no{\bf $3.2.$ An approach of simply objective decision}

\vskip 3mm

\no This approach is established under the following
considerations for Problems $1$ and $2$. \vskip 3mm

\no{\bf Consideration $1$} {\it In these evaluation items
$A_1,A_2,\cdots,A_m$ of a tendering $\widetilde{A}$, seek the
optimum of one evaluation item. For example, seek the lowest
bidding price in a construction tendering for a simply building or
seek the optimum of design scheme in a design project tendering,
etc.. }

\no{\bf Consideration $2$} {\it The value of these evaluation
items $A_1,A_2,\cdots,A_m$ is comparable which enables us to
measure each of them by a unify unit and to construct various
weighted functions on them. For example, the 100 marks and the
lowest evaluated price method widely used in China are used under
this consideration. }

\vskip 3mm

\no $3.2.1.$ The optimum of one objective

\vskip 2mm

\no Assume the optimal objective being $A_1$ in a tendering
$\widetilde{A}=\bigcup\limits_{i=1}^mA_i$. We need to determine
the acceptable basic criterions for all other items
$A_2,A_3,\cdots,A_k$, then arrange
$R_1(A_1),R_2(A_1),\cdots,R_l(A_1)$ among these acceptable bids
$R_1,R_2,\cdots,R_l$ for items $A_2,A_3,\cdots,A_k$ in $R_i, 1\leq
i\leq k$. For example, evaluating these items $A_2,A_3,\cdots,A_k$
by qualification or by weighted function on $A_2,A_3,\cdots,A_k$
up to these criterions, then arrange these acceptable bids
$R_1,R_2,\cdots,R_l$ under their response to $A_1$ and the order
of $R_i(\widetilde{A})$, $R_i(\widetilde{A})$ if $R_i(A_1)\approx
R_j(A_1)$. According to Theorem $3.1$, we get the following
result.

\vskip 4mm

\no{\bf Theorem $3.2$} \ {\it The approach of one optimal
objective can get an ordered sequence of bids for a tendering
$\widetilde{A}$.}

\vskip 3mm

\no{\bf Example $3.2$}\ {\it The optimum of design scheme is the
objective in a design project tendering $\widetilde{A}$ which is
divided into $5$ ranks $A,B,C,D,E$ and other evaluation items such
as human resources, design period and bidding price by a
qualifiable approach if the bidding price is in the interval of
the service fee norm of China. The final order of bids is
determined by the order of design schemes with qualifiable human
resources, design period and bidding price and applying the ruler
of lower price first for two equal design scheme in order.

There were $8$ bidders in this tendering. Their bidding prices are
in table $3$.

 \vskip 3mm

\begin{center}
\begin{tabular}{|c|c|c|c|c|c|c|c|c|}\hline
\ bidder \ & $R_1$ &  $R_2$ & $R_3$ & $R_4$ & $R_5$ & $R_6$ & $R_7$ & $R_8$  \\
\hline
 \ bidding \ price \
 & $251$ &  $304$ & $268$ & $265$ & $272$ & $283$ & $278$ & $296$\\ \hline
\end{tabular}
\end{center}
\vskip 2mm

\c{\rm table $3$}\vskip 2mm

After evaluation for these human resources, design period and
bidding price, $4$ bidders are qualifiable unless the bidder $R_5$
in  human resources. The evaluation result for bidding design
schemes is in table $4$. \vskip 3mm

\vskip 2mm
\begin{center}
\begin{tabular}{|c|c|c|c|c|c|}\hline
\ rank \ & $A$ &  $B$ & $C$ & $D$ & $E$  \\ \hline
 \ design \ scheme \
 & $R_3$¡¢$R_6$ &  $R_1$ & $R_2$¡¢$R_8$ & $R_7$ &  $R_4$ \\ \hline
\end{tabular}
\end{center}
\vskip 2mm

\c{\rm table $4$}

\vskip 2mm

Therefore, the ordered sequence for bids is}

$$R_3(\widetilde{A)\succ R_6(\widetilde{A})}\succ R_1(\widetilde{A})\succ
R_8(\widetilde{A})\succ R_2(\widetilde{A})\succ
R_7(\widetilde{A})\succ R_4(\widetilde{A}).$$

\vskip 3mm

\no{\bf Example $3.3$} \ {\it The optimum objective in a tendering
$\widetilde{A}$ for a construction of a dwelling house is the
bidding price $A_1$. All other evaluation items, such as
qualifications, management persons and equipments is evaluated by
a qualifiable approach.

There were $7$ bidders $R_i,1\leq i\le 7$ in this tendering. The
evaluation of price is by a weighted function approach, i.e.,
determine the standard price $S$ first, then calculate the mark
$N$ of each bidder by the following formulae}

$$S=\frac{(\sum\limits_{i=1}^7A_i-\max\{R_i(A_1)|1\leq i\leq 7\}
-\min\{R_i(A_1)|1\leq i\leq 7\}}{5},$$

$$N_i=100-t\times|\frac{R_i(A_1)-S}{S}|\times 100, \ \ 1\leq i\leq 7,$$

\no {\it where, if $R_i(A_1)-S > 0$ then $t=6$ and if $R_i(A_1)-S
< 0$ then $t=3$.

After evaluation, all bidders are qualifiable in qualifications,
management persons and equipments. Their bidding prices are in
table $5$.

\vskip 3mm
\begin{center}
\begin{tabular}{|c|c|c|c|c|c|c|c|}\hline
\ bidder \ & $R_1$ &  $R_2$ & $R_3$ & $R_4$ & $R_5$ & $R_6$ & $R_7$   \\
\hline
 \ $A_1$ \
 & $3518$ &  $3448$ & $3682$ & $3652$ & $3490$ & $3731$ & $3436$ \\ \hline
\end{tabular}
\end{center}
\vskip 2mm

\c{\rm table $5$}

\vskip 2mm

According to these formulae, we get that $S=3558$ and the mark of
each bidder as those shown in table $6$.

\vskip 3mm
\begin{center}
\begin{tabular}{|c|c|c|c|c|c|c|c|}\hline
\ bidder \ & $R_1$ &  $R_2$ & $R_3$ & $R_4$ & $R_5$ & $R_6$ & $R_7$   \\
\hline
 \ mark \
 & $96.70$ &  $91.27$ & $79.12$ & $84.16$ & $94.27$ & $73.84$ & $89.68$ \\ \hline
\end{tabular}
\end{center}
\vskip 2mm

\c{\rm table $6$}\vskip 2mm

Therefore, the ordered sequence of bids is}

$$R_1(\widetilde{A})\succ R_5(\widetilde{A})\succ R_2(\widetilde{A})\succ
R_7(\widetilde{A})\succ R_4(\widetilde{A})\succ
R_3(\widetilde{A})\succ R_6(\widetilde{A}).$$

\vskip 3mm

\no $3.2.2.$ The pseudo-optimum of multiple objectives

\vskip 2mm

\no This approach assumes that there is a unifying unit between
these evaluation items $A_1,A_2,\cdots,A_m$ in an interval
$[a,b]$. Whence it can be transformed into case $3.2.1$ and sought
the optimum of one objective. Not loss of generality, we assume
the unifying unit is $\varpi$ and

$$\varpi(A_i)=f_i(\varpi), \ \ 1\leq i\leq m,$$

\no where $f_i$ denotes the functional relation of the metric
$\varpi(A_i)$ with unit $\varpi$. Now the objective of tendering
turns to a programming of one objective

$$\max\limits_{\varpi}F(f_1(\varpi),f_2(\varpi),\cdots,f_m(\varpi)) \ {\rm or} \
\min\limits_{\varpi}F(f_1(\varpi),f_2(\varpi),\cdots,f_m(\varpi)),$$

\no where $F$ denotes the functional relation of the tendering
$\widetilde{A}$ with these evaluation items $A_1,A_2,\cdots,A_m$,
which can be a weighted function, such as a linear function

$$F(f_1(\varpi),f_2(\varpi),\cdots,f_m(\varpi))=
\sum\limits_{i=1}^mf_i(\varpi)$$

\no or an ordered sequence. According to Theorem $3.2$, we know
the following result.

\vskip 4mm

\no{\bf Theorem $3.3$} \ {\it If the function $F$ of a tendering
$\widetilde{A}$ only has one maximum value in $[a,b]$, then there
exists an ordered sequence for these bids
$R_i(\widetilde{A}),1\leq i\leq k$ after determined how to arrange
$R_i(\widetilde{A})$ and $R_j(\widetilde{A})$ when
$F(R_i(\widetilde{A}))=F(R_j(\widetilde{A})),i\not=j$.}

\vskip 3mm

The 100 marks and the lowest evaluated price method widely used in
China both are applications of this approach. In the 100 marks,
the weight function is a linear function

$$F(f_1(\varpi),f_2(\varpi),\cdots,f_m(\varpi))=
\sum\limits_{i=1}^mf_i(\varpi)£¬$$

\no with $0\leq F(f_1(\varpi),f_2(\varpi),\cdots,f_m(\varpi))\leq
100, f_i\geq 0, 1\leq i\leq m$. In the lowest evaluated price
method, each difference of an evaluation item $A_i, 2\leq i\leq m$
is changed to the bidding price $\varpi(A_1)$, i.e.,

$$f_i=(R(A_i)-S(A_i))\varpi(A_1), \ 1\leq i\leq m,$$

\no where $S(A_i)$ is the standard line for $A_i$, $\varpi(A_i)$
is one unit difference of $A_i$ in terms of $A_1$. The weighted
function of the lowest evaluated price method is

$$F(\varpi(A_1),f_2(\varpi(A_1)),\cdots,f_m(\varpi(A_1)))=(1+\sum\limits_{j=2}^m(R(A_i)-S(A_i)))\varpi(A_1).$$

For example, we can fix one unit difference of a technological
parameter $15$, i.e., $\varpi(A_1)=15$ ten thousand dollars in
terms of the bidding price.

\vskip 5mm

\no{\bf \S $4.$ Weighted functions and their construction}\vskip
3mm

\no We discuss weighted functions on the evaluation items or
indexes in this section. First, we give a formal definition for
weighted functions.

\vskip 4mm

\no{\bf Definition $4.1$} \ {\it For a tendering
$\widetilde{A}=\bigcup\limits_{i=1}^mA_i$, where
$A_i=\{a_{i1},a_{i2},\cdots,a_{in}\}, 1\leq i\leq m$ with $k$
bidders $R_1,R_2,\cdots,R_k$, if there is a continuous function
$\omega:\widetilde{A}\rightarrow [a,b]\subset(-\infty,+\infty)$ or
$\omega:A_i\rightarrow [a,b]\subset(-\infty,+\infty),1\leq i\leq
m$ such that for any integers $l,s,1\leq l,s\leq k$,
$R_l(\omega(\widetilde{A}))>R_s(\omega(\widetilde{A}))$ or
$R_l(\omega(\widetilde{A}))=R_s(\omega(\widetilde{A}))$ as
$R_l(\widetilde{A})\succ R_s(\widetilde{A})$ or
$R_l(\widetilde{A})\approx R_s(\widetilde{A})$ and
$R_l(\omega(A_i)>R_s(\omega(A_i))$ or
$R_l(\omega(A_i))=R_s(\omega(A_i))$ as $R_l(A_i)\succ R_s(A_i)$ or
$R_l(A_i)\approx R_s(A_i), 1\leq i\leq m$, then $\omega$ is called
a weighted function for the tendering $\widetilde{A}$ or the
evaluation items $A_i,1\leq i\leq m$.}

\vskip 3mm

According to the decision theory of multiple objectives($[3]$),
the weighted function $\omega(A_i)$ must exists for any integer
$i, 1\leq i\leq m$. but generally, the weight function
$\omega(\widetilde{A})$ does not exist if the values of these
evaluation items $A_1,A_2,\cdots,A_m$ can not compare. There are
two choice for the weighted function $\omega(A_i)$.\vskip 3mm

\no{\bf Choice $1$} \ {\it the monotone functions in the interval
$[a,b]$, such as the linear functions.}\vskip 2mm

\no{\bf Choice $2$} \ {\it The continuous functions only with one
maximum value in the interval $[a,b]$, such as
$\omega(A_i)=-2x^2+6x+12$ or}

\[
\omega(A_i)=\left\{\begin{array}{lr}
x,& {\rm if}\quad 0\leq x\leq 2,\\
-x+4,& {\rm if}\quad x\geq 4.
\end{array}
\right.
\]

As examples of concrete weighted functions $\omega$, we discuss
the tendering of civil engineering constructions.

\vskip 4mm

\no{\bf $4.1.$ The weighted function for the bidding price}

\vskip 3mm

\no Let $A_1$ be the bidding price. We often encounter the
following weighted function $\omega(A_1)$ in practice.

$$\omega(R_i(A_1))=-\varsigma\times\frac{R_i(A_1)-S}{S}+\zeta$$

\no where,

$$
S=\frac{R_1(A_1)+R_2(A_1)+\cdots +R_k(A_1)}{k}
$$

\no or

\[
S=\left\{
\begin{array}{ll}
\frac{R_1(A_1)+R_2(A_1)+\cdots +R_k(A_1)-M-N}{k-2},& k\geq 5,\\
\frac{R_1(A_1)+R_2(A_1)+\cdots +R_k(A_1)}{k}, & 3\leq k\leq 4 \\
\end{array}
\right.
\]

\no or

$$S = T\times A\% +\frac{R_1(A_1)+R_2(A_1)+\cdots+R_k(A_1)}{k}\times
(1-A\%).$$\vskip 2mm

\no Where $T$,$A$\%,$k$, $M$ and $N$ are the pre-price of the
tender, the percentage of $T$ in $S$, the number of bidders and
the maximum and minimum bidding price, respectively,
$R_i(A_1),i=1,2,\cdots,k$ denote the bidding prices and
$\varsigma$, $\zeta$ are both constants.

There is a postulate in these weighted functions, i.e., each
bidding price is random and accord with the normal distribution.
Then the best excepted value of this civil engineering is the
arithmetic mean of these bidding prices. However, each bidding
price is not random in fact. It reflects the bidder's expected
value and subjectivity in a tendering. We can not apply any
definite mathematics to fix its real value. Therefore, this
formula for a weighted function can be only seen as a game, not a
scientific decision.

By the view of scientific decision, we can apply weighted
functions according to the expected value and its cost in the
market, such as

\vskip 3mm

($1$) \ {\it the linear function

$$\omega(R_i(A_1))=-p\times\frac{R_i(A_1)-N}{M-N}+q$$

\no in the interval $[N,M]$, where $M,N$ are the maximum and
minimum bidding prices $p$ is the deduction constant and $q$ is a
constant such that $R_i(\omega(A_1))\geq 0, 1\leq i\leq k$. The
objective of this approach is seek a lower bidding price. }\vskip
3mm

($2$) \ {\it non-linear functions in the interval $[N,M]$, such as

$$\omega(R_i(A_1))=-p\times\frac{R_i(A_1)-\frac{T+\sum\limits_{j=1}^kR_i(A_1)}{k+1}}+q,$$

$$\omega(R_i(A_1))=-p\times\frac{R_i(A_1)-\sqrt[k+1]{R_1(A_1)R_2(A_1)\cdots
R_k(A_1)T}} {\sqrt[k+1]{R_1(A_1)R_2(A_1)\cdots R_k(A_1)T}}+q$$

\no or

$$\omega(R_i(A_1))=-p\times\frac{R_i(A_1)-
\sqrt{\frac{R_1^2(A_1)+R_2^2(A_1)+\cdots+R_k^2(A_1)+T^2}{k+1}}}
{\sqrt{\frac{R_1^2(A_1)+R_2^2(A_1)+\cdots+R_k^2(A_1)+T^2}{k+1}}}+q$$

\no etc..} If we wish to analog a curve for these bidding prices
and choose a point on this curve as $\omega(R_i(A_1))$, we can
apply the value of a polynomial of degree $k+1$

$$f(x)=a_{k+1}x^{k+1}+a_kx^k+\cdots+a_1x+a_0$$

\no by the undetermined coefficient method. Arrange the bidding
prices and pre-price of the tender from lower to higher. Not loss
of generality, let it be $R_{j_1}(A_1)\succ
R_(j_2)(A_1)\succ\cdots\succ T\succ\cdots\succ R_{j_k}(A_1)$.
Choose $k+2$ constants $c_1> c_2> \cdots> c_{k+1}> 0$, for
instance $k+1>k>\cdots> 1> 0$. Solving the equation system

\begin{eqnarray*}
& \ & R_{j_1}(A_1)=a_{k+1}c_{1}^{k+1}+a_kc_1^k+\cdots+a_1c_1+a_0\\
& \ & R_{j_2}(A_1)=a_{k+1}c_{2}^{k+1}+a_kc_2^k+\cdots+a_1c_2+a_0\\
& \ & \cdots\cdots\cdots\cdots\cdots\cdots\cdots\cdots\cdots\\
& \ &
R_{j_{k-1}}(A_1)=a_{k+1}c_{k}^{k+1}+a_kc_{k}^k+\cdots+a_1c_k+a_0\\
& \ & R_{j_k}(A_1)=a_0
\end{eqnarray*}

\no we get a polynomial $f(x)$ of degree $k+1$. The bidding price
has an acceptable difference in practice. Whence, we also need to
provide a bound for the difference which does not affect the
ordered sequence of bids.

\vskip 4mm

\no{\bf $4.2.$ The weighted function for the programming}

\vskip 3mm

\no Let $A_2$ be the evaluation item of programming with
evaluation indexes $\{a_{21},a_{22},$ $\cdots,a_{2n_2}\}$. It is
difficult to evaluating a programming in quantify, which is not
only for the tender, but also for the evaluation specialists. In
general, any two indexes of $A_2$ are not comparable. Whence it is
not scientific assigning numbers for each index since we can not
explain why the mark of a programming is $96$ but another is $88$.
This means that it should qualitatively evaluate a programming or
a quantify after a qualitatively evaluation. Its weight function
$\omega(R_i(A_2)),1\leq i\leq k$ can be chosen as a linear
function

$$\omega(R_i(A_2))=\omega(R_i(a_{21}))+\omega(R_i(a_{22}))+\cdots+\omega(R_i(a_{2n_2})).$$

For example, there are $4$ evaluation indexes for the programming,
and each with $A,B,C,D$ ranks in a tendering. The corespondent
mark for each rank is in table $7$.

\vskip 3mm
\begin{center}
\begin{tabular}{|c|c|c|c|c|}\hline
\ index \ & $a_{21}$ &  $a_{22}$ & $a_{23}$ & $a_{24}$     \\
\hline
 \ $A$ \ & $ \ 4 \ $ &  $ \ 2 \ $ & $ \ 2 \ $ & $ \ 1 \ $ \\ \hline
 \ $B$ \ & $3$ &  $1.5$ & $1.5$ & $0.8$ \\ \hline
 \ $C$ \ & $2$ &  $1$ & $1$ & $0.5$ \\ \hline
 \ $D$ \ & $1$ &  $0.5$ & $0.5$ & $0.3$ \\ \hline
\end{tabular}
\end{center}
\vskip 3mm

\c{\rm table $7$}\vskip 2mm

If the evaluation results for a bidding programming $R_i, 1\leq
i\leq 4$ are $\omega(R_i(a_{21}))=A$, $\omega(R_i(a_{22}))=B$,
$\omega(R_i(a_{23}))$ $=B$ and $\omega(R_i(a_{24}))=A$, then the
mark of this programming is

\begin{eqnarray*}
R_i(\omega(A_2))&=& R_i(\omega(a_{21}))+R_i(\omega(a_{22}))+
R_i(\omega(a_{23}))+R_i(\omega(a_{24}))\\
&=& 4+3+1.5+1=9.5.
\end{eqnarray*}

By the approach in Section $3$, we can alphabetically or graphicly
arrange the order of these programming if we can determine the
rank of each programming. Certainly, we need the order of these
indexes for a programming first. The index order for programming
is different for different constructions tendering.

\vskip 5mm

\no{\bf \S $5.$ Further discussions}\vskip 4mm

\no{\bf $5.1$} \ Let $\widetilde{A}=\bigcup\limits_{i=1}^mA_i$ be
a Smarandache multi-space with an operation set
$O(\widetilde{A})=\{\circ_i; 1\leq i\leq m\}$. If there is a
mapping $\Theta$ on $\widetilde{A}$ such that
$\Theta(\widetilde{A})$ is also a Smarandache multi-space, then
$(\widetilde{A},\Theta)$ is called a {\it pseudo-multi-space}.
Today, nearly all geometries, such as the {\it Riemann} geometry,
{\it Finsler} geometry and these pseudo-manifold geometries  are
particular cases of pseudo-multi-geometries.

For applying Smarandache multi-spaces to an evaluation system,
choose $\Theta(\widetilde{A})$ being an order set. Then Theorem
$3.1$ only asserts that any subset of $\Theta(\widetilde{A})$ is
an order set, which enables us to find the ordered sequence for
all bids in a tendering. Particularly, if $\Theta(\widetilde{A})$
is continuous and $\Theta(\widetilde{A})\subseteq
[-\infty,+\infty]$, then $\Theta$ is a weighted function on
$\widetilde{A}$ widely applied in the evaluation of bids in China.
By a mathematical view, many problems on $(\widetilde{A},\Theta)$
is valuable to research. Some open problems are presented in the
following.

\vskip 3mm

\no{\bf Problem $5.1$} \ {\it Characterize these
pseudo-multi-spaces $(\widetilde{A},\Theta)$, particularly, for
these cases of
$\Theta(\widetilde{A})=\bigcup\limits_{i=1}^n[a_i,b_i]$,
$\Theta(\widetilde{A})=\bigcup\limits_{i=1}^n(G_i,\circ_i)$ and
$\Theta(\widetilde{A})=\bigcup\limits_{i=1}^n(R;+_i,\circ_i)$ with
$(G_i,\circ_i)$ and $(R;+_i,\circ_i)$ being a finite group or a
ring for $1\leq i\leq n$.}

\vskip 3mm

\no{\bf Problem $5.2$}\ {\it Let $\Theta(\widetilde{A})$ be a
group, a ring or a filed. Can we find an ordered sequence for a
finite subset of $\widetilde{A}$?}

\vskip 3mm

\no{\bf Problem $5.3$} \ {\it Let $\Theta(\widetilde{A})$ be $n$
lines or $n$ planes in an Euclid space ${\bf R}^n$. Characterize
these pseudo-multi-spaces $(\widetilde{A},\Theta)$. Can we find an
arrangement for a finite subset of $\widetilde{A}$?}

\vskip 3mm

\no{\bf $5.2$} \ The evaluation approach in this paper can be also
applied to evaluate any multiple objectives, such as the
evaluation of a scientific project, a personal management system,
an investment of a project, $\cdots$, etc..

\vskip 10mm

\no{\bf References}\vskip 4mm

\re{[1]}G.Chartrand and L.Lesniak,{\it Graphs \& Digraphs},
Wadsworth, Inc., California, 1986.

\re{[2]}T.Chen, {\it Decision Approaches in Multiple
Objectives}(in Chinese), in {\it The Handbook of Modern
Engineering Mathematics}, Vol.IV, Part 77,1357-1410. Central China
Engineering College Press, 1987.

\re{[3]}P.C.Fishburn, {\it Utility Theory for Decision Making},
New York, Wiley, 1970.

\re{[4]}D.L.Lu, X.S.Zhang and Y.Y.Mi, An offer model for civil
engineering construction, {\it Chinese OR Transaction}, Vol.5,
No.4(2001)41-52.

\re{[5]}L.F.Mao, {\it On Automorphisms groups of Maps, Surfaces
and Smarandache geometries}, {\it Sientia Magna}, Vol.$1$(2005),
No.$2$, 55-73.

\re{[6]}L.F.Mao, {\it Automorphism Groups of Maps, Surfaces and
Smarandache Geometries}, American Research Press, 2005.

\re{[7]}L.F.Mao, {\it Smarandache multi-space theory}, Hexis,
Phoenix, AZ£¬2006.

\re{[8]}L.F.Mao, {\it Chinese Construction Project Bidding
Technique \& Cases Analyzing--Smarandache Multi-Space Model of
Bidding},Xiquan Publishing House (Chinese Branch), America, 2006.

\re{[9]}F.Smarandache, Mixed noneuclidean geometries, {\it eprint
arXiv: math/0010119}, 10/2000.

\re{[10]}F.Smarandache, {\it A Unifying Field in Logics.
Neutrosopy: Neturosophic Probability, Set, and Logic}, American
research Press, Rehoboth, 1999.

\re{[11]}F.Smarandache, Neutrosophy, a new Branch of Philosophy,
{\it Multi-Valued Logic}, Vol.8, No.3(2002)(special issue on
Neutrosophy and Neutrosophic Logic), 297-384.

\re{[12]}F.Smarandache, A Unifying Field in Logic: Neutrosophic
Field, {\it Multi-Valued Logic}, Vol.8, No.3(2002)(special issue
on Neutrosophy and Neutrosophic Logic), 385-438.

\end{document}